\documentstyle[11pt,epsf]{article}
\title{Interlacing properties and the Schur-Szeg\H{o} composition}
\author{Vladimir P. Kostov
\footnote{Research partially supported by research project 20682 
for cooperation between CNRS and FAPESP 
``Zeros of algebraic polynomials''}} 
\date{}
\newtheorem{tm}{Theorem}
\newtheorem{defi}[tm]{Definition}
\newtheorem{rem}[tm]{Remark}

\newtheorem{lm}[tm]{Lemma}

\newtheorem{cor}[tm]{Corollary}
\newtheorem{prop}[tm]{Proposition}
\newtheorem{nota}[tm]{Notation}
\begin{document}
\maketitle 

\begin{abstract}
Each degree $n$ polynomial in one variable of the form 
$(x+1)(x^{n-1}+c_1x^{n-2}+\cdots +c_{n-1})$ is representable in a unique way 
as a Schur-Szeg\H{o} composition of $n-1$ polynomials of the form 
$(x+1)^{n-1}(x+a_i)$, see \cite{Ko1}, \cite{AlKo} and \cite{Ko2}. Set 
$\sigma _j:=\sum _{1\leq i_1<\cdots <i_j\leq n-1}a_{i_1}\cdots a_{i_j}$. 
The eigenvalues of the affine mapping 
$(c_1,\ldots ,c_{n-1})\mapsto (\sigma _1,\ldots ,\sigma _{n-1})$ 
are positive rational numbers and its eigenvectors are defined by 
hyperbolic polynomials (i.e. with real roots only). In the present paper 
we prove interlacing properties of the roots of these polynomials. 

{\bf Key words:} Schur-Szeg\H{o} composition; composition factor

{\bf AMS classification:} 12D10
\end{abstract}

\section{The Schur-Szeg\H{o} composition and the mapping $\Phi$}

The present paper deals with the {\em composition of Schur-Szeg\H{o} (CSS)} 
of degree $n$ polynomials in one 
variable. For the polynomials $P:=\sum _{j=0}^na_jx^j$ and 
$Q:=\sum _{j=0}^nb_jx^j$ their CSS is defined by 
the formula $P\ast Q:=\sum _{j=0}^na_jb_jx^j/C_n^j$ ($C_n^j=n!/j!(n-j)!$, 
$a_i, b_i\in {\bf C}$). 
In the monographies \cite{Pr} and 
\cite{RS} one can find properties and 
applications of the CSS. The CSS is associative, commutative and polylinear. 
It can be defined for more than 
two polynomials by the formula  

\[ P_1\ast \cdots \ast P_s=\sum _{j=0}^na_j^1\cdots a_j^sx^j/(C_n^j)^{s-1}~~
{\rm where}~~P_i=\sum _{j=0}^na_j^ix^j~.\] 

The role of the unity in the CSS is played by the polynomial $(x+1)^n$ 
in the sense that for every degree $n$ polynomial $P$ one has 
$P\ast (x+1)^n=P$. 

\begin{defi}
{\rm A real polynomial is {\em hyperbolic} (resp. {\em strictly hyperbolic}) 
if it has only real (resp. only real and distinct) roots.}
\end{defi}

In papers \cite{KoSh}, \cite{Ko1} and \cite{Ko4} the question is considered 
how many of the roots of the CSS of two hyperbolic or only real 
polynomials are real 
negative, zero or positive. In the case of hyperbolic polynomials the 
exhaustive answer is given in \cite{Ko1}, while \cite{Ko4} contains 
sufficient conditions for the realizability of certain cases defined by the 
number of negative, positive and complex roots of $P$, $Q$ and $P\ast Q$. 
In the case of two hyperbolic polynomials $P$ and $Q$ one of which has only 
negative roots, the multiplicity vector of the roots of $P\ast Q$ is completely 
defined by the ones of $P$ and $Q$, see \cite{KoSh}. 

Call {\em composition factor} any polynomial of the 
form $K_a:=(x+1)^{n-1}(x+a)$. In paper \cite{Ko1} it is announced and in 
paper \cite{AlKo} it is proved that any monic 
degree $n$ complex 
polynomial $P$ such that $P(-1)=0$ is representable in the form 

\begin{equation}\label{factoriz}
P=K_{a_1}\ast \cdots \ast K_{a_{n-1}}
\end{equation}
where the numbers $a_i\in {\bf C}$ are unique up to permutation. To extend 
this presentation to the case when $P$ is not necessarily monic and, more 
generally, to the case when deg$P\leq n$, one has to admit the presence of 
composition factors  
$K_{\infty}:=(x+1)^{n-1}$ and of a constant factor in the right-hand side of 
(\ref{factoriz}). When 
the polynomial $P$ is real, then part of the numbers $a_i$ are real while the 
rest form complex conjugate couples. (Otherwise conjugation of 
(\ref{factoriz}) defines a new $(n-1)$-tuple of numbers $a_i$.) 

Denote by $\sigma _j$ the $j$th elementary symmetric polynomial of the numbers 
$a_i$, i.e. 
$\sigma _j:=\sum _{1\leq i_1<\cdots <i_j\leq n-1}a_{i_1}\cdots a_{i_j}$. Set 
$P:=(x+1)(x^{n-1}+c_1x^{n-2}+\cdots +c_{n-1})$. Consider the mapping  

\[ \Phi ~:~(c_1,\ldots ,c_{n-1})~
\mapsto ~(\sigma _1,\ldots ,\sigma _{n-1})~.\]  

In paper \cite{Ko2} it is shown that this mapping is affine, that its 
eigenvalues are the positive rational numbers listed in the first line of the 
following displayed formula, and the corresponding eigenvectors are defined 
by the polynomials of its second line:

\[ \begin{array}{cccccc}\lambda _1=1&\lambda _2=\frac{n}{(n-1)}&
\lambda _3=\frac{n^2}{(n-1)(n-2)}&
\lambda _4=\frac{n^3}{(n-1)(n-2)(n-3)}&\ldots &
\lambda _{n-1}=\frac{n^{n-2}}{(n-1)!}\\ \\  
(x+1)^{n-1}&x(x+1)^{n-2}&x(x+1)^{n-3}Q_1(x)&x(x+1)^{n-4}Q_2(x)&\ldots &
x(x+1)Q_{n-3}(x)\end{array}\]  

\begin{rem}\label{remlin}
{\rm 1) The polynomial $Q_j$ is degree $j$, monic, strictly 
hyperbolic, with all roots positive and {\em self-reciprocal}, i.e. 
$x^jQ_j(1/x)=\pm Q_j(x)$. In particular, $Q_1=x-1$. Hence if $x_0$ is root of 
$Q_j$, then $1/x_0$ is also such a root. If $j$ is odd, then $Q_j(1)=0$.

2) In the above list the ``eigenpolynomials'' are all of degree $n-1$ 
and divisible by $(x+1)$. If one considers $\Phi$ as a linear (not affine) 
mapping $(c_0,\ldots ,c_{n-1})\mapsto (\sigma _0,\ldots ,\sigma _{n-1})$, 
i.e. when $P=(x+1)(c_0x^{n-1}+c_1x_{n-2}+\cdots +c_{n-1})$, then one sets  
$\sigma _0=c_0$ and has to add a second eigenvalue $1$ to which there 
corresponds the eigenpolynomial $(x+1)^n$. To the eigenvalue $1$ there 
correspond two Jordan blocks of size $1$. 

3) It is shown in \cite{KoMaSh} that for $j$ fixed and when 
$n\rightarrow \infty$, the polynomial $xQ_j(-x)$ tends to the 
{\em Narayana polynomial} $\sum _{i=1}^nN_{n,i}x^i$ where $N_{n,i}$ are the 
{\em Narayana numbers} $C_n^iC_n^{i-1}/n$.}
\end{rem}

Denote by $0<x_1<\cdots <x_j$ the roots of $Q_j$, by 
$0<y_1<\cdots <y_{j+1}$ the ones of $Q_{j+1}$ and by $0<z_1<\cdots <z_{j+2}$ 
the ones of $Q_{j+2}$. In the present paper we prove the following two 
theorems. In their proofs we use some known results about 
hyperbolic polynomials and the mapping $\Phi$, see Section~\ref{prelim}. 

\begin{tm}\label{tm1}
(First interlacing property) The roots of $Q_j$ and $Q_{j+1}$ interlace 
($j=1,\ldots ,n-4$), i.e. 
$y_1<x_1<y_2<x_2<\cdots <y_j<x_j<y_{j+1}$.
\end{tm}

The theorem is proved in Section~\ref{prtm1}.

\begin{tm}\label{tm2} 
(Second interlacing property) For $i=1,\ldots ,[j/2]$ one has 
$x_i\in (z_i,z_{i+1})$ and $x_{j+1-i}\in (z_{j+2-i},z_{j+3-i})$. 
If $j$ is odd, then 
$x_{(j+1)/2}=z_{(j+3)/2}=1$.
\end{tm}

The theorem is proved in Section~\ref{prtm2}. An analog of Theorem~\ref{tm2} 
holds for the roots of Gegenbauer polynomials and their derivatives 
and for the roots of Narayana 
polynomials, see Section~\ref{Geg}.

\section{Preliminaries\protect\label{prelim}}

In the present section we recall some classical properties of hyperbolic 
polynomials and some properties of the mapping $\Phi$. The following theorem 
is a well-known result which is used in the proofs (see Section~5 in \cite{O}):

\begin{tm}\label{preplit}
Suppose that $P$ and $Q$ are two polynomials with no root in common. Then 
they are hyperbolic and their roots interlace if and only if for any 
$(\theta ,\mu )\in ({\bf R}^2\backslash \{ (0,0)\} )$ the polynomial 
$\theta P+\mu Q$ has only real roots.
\end{tm}

Next, we recall some nontrivial properties of the mapping $\Phi$:

\begin{prop}\label{multiplicity}
If the degree $n$ polynomials $P$ and $Q$ have nonzero roots 
$x_P$, $x_Q$ of multiplicities $m_P$, $m_Q$ such that $m_P+m_Q\geq n$, 
then $-x_Px_Q$ is a root of $P\ast Q$ of multiplicity $m_P+m_Q-n$.
\end{prop}

This is Proposition~1.4 in \cite{KoSh} 
(in \cite{KoSh} the condition the roots to be nonzero is omitted which is 
not correct).

\begin{rem}\label{remmult}
{\rm In the right-hand side of (\ref{factoriz}) there are exactly $l$ 
composition factors with $a_i\neq -1$ if and only if 
the multiplicity of $(-1)$ 
as a root of $P$ equals $n-l$. This follows from 
Proposition~\ref{multiplicity} applied to the right-hand side $l-1$ times, 
the roles of $x_P$, $x_Q$ being played by the roots equal to $(-1)$.}
\end{rem}

\begin{nota}
{\rm For $i=0, 1,\ldots ,n-1$ we set $b_i:=-i/(n-i)$. Set $b_n:=-\infty$.}
\end{nota}

\begin{prop}\label{propmult}
1) Suppose that $P=x^mR$, deg$R\leq n-m$.
Then $m$ of the numbers $a_i$ equal (up to permutation) 
$b_0$, $b_1$, $\ldots$, $b_{m-1}$. 

2) Suppose that deg$P=n-m$. Then $m$ of the numbers $a_i$ equal 
(up to permutation) $b_n$, $b_{n-1}$, $\ldots$, $b_{n-m+1}$.
\end{prop}

The proposition coincides with part 3) of Remark~2 in \cite{Ko2}.

\begin{prop}\label{compfactors}
For $a<0$ there is exactly one change of sign 
in the sequence of coefficients of the 
polynomial $K_a$. For $a>0$ there is no change of sign in it. 
\end{prop}

Indeed, one has 

\[ K_a=(x+1)^{n-1}(x+a)=\sum _{j=0}^n(aC_{n-1}^j+C_{n-1}^{j-1})x^j=
x^n+\sum  _{j=0}^{n-1}C_{n-1}^j\left( a+\frac{j}{n-j}\right) x^j \] 
and the signs of the sequence of coefficients are the same as the ones of the 
sequence $a$, $a+1/(n-1)$, $a+2/(n-2)$, $\ldots$, $a+n-1$, $1$.~~~~~$\Box$

\begin{prop}\label{posroots}
If the hyperbolic polynomial $P$ has $l$ positive roots, then at least $l$ 
of the numbers $a_i$ are negative and distinct.
\end{prop}

The proposition follows from the final remark in \cite{Ko3}.

\begin{prop}\label{preserverecipr}
The mapping $\Phi$ preserves self-reciprocity.
\end{prop}

This follows from the fact that the matrix of $\Phi$ (considered as a linear, 
not affine mapping, see part 2) of Remark~\ref{remlin}) is centre-symmetric, 
see Proposition~10 in \cite{Ko2}.

\section{Proof of Theorem~\protect\ref{tm1}\protect\label{prtm1}}

$1^0$. Consider two polynomials of the form $U:=x(x+1)^{n-2-j}T_j(x)$ and 
$V:=x(x+1)^{n-3-j}T_{j+1}(x)$ where the polynomials $T_i$ are monic, 
hyperbolic and self-reciprocal, deg$T_i=i$. Assume also that the roots of 
$T_j$ and $T_{j+1}$ are all distinct, positive and that they interlace 
(in the same sense as in the claim of Theorem~\ref{tm1} about the roots 
of $Q_j$ and $Q_{j+1}$). 

\begin{prop}\label{propinterlace}
Under these assumptions, for every 
$(\theta ,\mu )\in ({\bf R}^2\backslash \{ (0,0)\} )$ 
the polynomial $T:=(x+1)\theta T_j+\mu T_{j+1}$ 
has all roots real and distinct at least $j$ of which are positive. 
The last root might be 
positive, negative or $0$ (in the last case the polynomial 
$\theta U+\mu V=x(x+1)^{n-3-j}T$ has a double 
root at $0$); it might equal $(-1)$ in which case the polynomial 
$\theta U+\mu V$ has an $(n-2-j)$-fold root at $(-1)$.
\end{prop}

{\em Proof:}\\ 

The cases $\theta =0\neq \mu$ and $\theta \neq 0=\mu$ are self-evident, 
so suppose that 
$\theta >0$, $\mu \neq 0$. 
The roots of the polynomials $(x+1)T_j$ and 
$T_{j+1}$ interlace. Hence in every interval between two consecutive roots of 
$T_{j+1}$ the polynomial $T$ 
takes at least one positive and at least one negative value. 

Indeed, it suffices to prove this for the interval  
between the two largest roots $0<a<b$ of $T_{j+1}$. Denote by $c$ the root 
of $T_j$ belonging to $(a,b)$. If $\mu >0$, then $T(b)>0$, $T(c)<0$. If 
$\mu <0$, then $T(c)>0$, $T(a)<0$.

Hence in all cases the polynomial $T$ has a root in $(a,b)$. In the same 
way one shows that it has a root between any two consecutive roots of 
$T_{j+1}$, i.e. it has at least $j$ positive distinct roots. 
For $(\theta ,\mu )=(0,1)$ 
its last root is positive, for $(\theta ,\mu )=(1,0)$ it equals $(-1)$. It is 
clear that by arguments of continuity there exists $(\theta ,\mu )\neq (0,0)$ 
for which the last root equals $0$.~~~~~$\Box$

\begin{prop}\label{hypint}
For $U$ and $V$ as above and for every 
$(\theta ,\mu )\in ({\bf R}^2\backslash \{ (0,0)\} )$ the 
polynomial $\Phi (\theta U+\mu V)$ is hyperbolic. It has a simple or double 
root at $0$, a root at $(-1)$ of multiplicity $(n-j-3)$ or $(n-j-2)$ and 
at least $j$ distinct positive roots.
\end{prop}

Indeed, by Proposition~\ref{propinterlace} the polynomial $\theta U+\mu V$ 
is hyperbolic and has at least $j$ distinct positive roots. 
By Remark~\ref{remmult} it has an $(n-j-3)$- 
or $(n-j-2)$-fold root at $(-1)$. 

By part 1) of Proposition~\ref{propmult} the polynomial 
$\Phi (\theta U+\mu V)$ has a root at $0$; by Proposition~\ref{posroots} 
it has $j$ or $j+1$ distinct positive roots. 
Hence in all cases it is hyperbolic. Indeed, it is of degree $n-1$, and 
except the positive distinct roots and the ones at $0$ and $(-1)$, there  
remains only one root not accounted for which  
is necessarily real.~~~~~$\Box$\\ 

$2^0$. The above proposition and Theorem~\ref{preplit} imply that the roots of 
the polynomials $\Phi (U)/x(x+1)^{n-3-j}$ and $\Phi (V)/x(x+1)^{n-3-j}$ 
are real and interlace.
Using the same type of arguments one sees that in the special cases 
$(\theta ,\mu )=(1,0)$ or $(0,1)$ the multiplicity of 
the root at $(-1)$ of $\Phi (\theta U+\mu V)$ equals respectively 
$(n-2-j)$ and $(n-3-j)$. In both cases the root at $0$ is simple. 
Hence there are respectively $j$ and $j+1$ distinct positive roots. 

$3^0$. This means that the polynomials $\Phi (U)$, $\Phi (V)$ satisfy the 
conditions of Proposition~\ref{hypint}. 
Hence the conclusion of the proposition 
applies also to the polynomial $\Phi ^2(\theta U+\mu V)$ and in the 
same way to $\Phi ^k(\theta U+\mu V)$ for $k\in {\bf N}$, i.e. we have 

\begin{cor}\label{corcor} 
For $k\in {\bf N}$, for $U$ and $V$ as above and for every 
$(\theta ,\mu )\in ({\bf R}^2\backslash \{ (0,0)\} )$ the 
polynomial $\Phi ^k(\theta U+\mu V)$ is hyperbolic. It has a simple or double 
root at $0$, a root at $(-1)$ of multiplicity $(n-j-3)$ or $(n-j-2)$ and 
at least $j$ distinct positive roots.
\end{cor}

\begin{prop}\label{zeroset}
All polynomials $\Phi ^k(U)$ and $\Phi ^k(V)$ are self-reciprocal. 
The set of roots of the polynomial $\Phi ^k(U)$ (resp. $\Phi ^k(V)$) 
tends for $k\rightarrow \infty$ to the one of 
$x(x+1)^{n-2-j}Q_j(x)$ (resp. the one of $x(x+1)^{n-3-j}Q_{j+1}(x)$).
\end{prop}

{\em Proof:}\\ 

Self-reciprocity of $\Phi ^k(U)$ and $\Phi ^k(V)$ follows from 
Proposition~\ref{preserverecipr}. Hence about 
half of the roots of the polynomials 
$\Phi ^k(U)/x(x+1)^{n-2-j}$ and $\Phi ^k(V)/x(x+1)^{n-3-j}$ 
belong to the interval $(0,1)$. 

Set $W_j:=x(x+1)^{n-2-j}Q_j$. The polynomials $W_i$, $i=0,\ldots ,j$, are a 
basis of the linear space of degree $\leq n-1$ polynomials divisible 
by $x(x+1)^{n-2-j}$. Set $U:=\sum _{i=0}^j\alpha _iW_i$, 
$V:=\sum _{i=0}^{j+1}\beta _iW_i$. One has $\alpha _j\neq 0\neq \beta _{j+1}$ 
because otherwise $U$ (resp. $V$) is divisible by 
$(x+1)^{n-1-j}$ (resp. by $(x+1)^{n-2-j}$). 

Hence $\Phi ^k(U)=\sum _{i=0}^j\alpha _i\lambda _i^kW_i$, 
$\Phi ^k(V)=\sum _{i=0}^{j+1}\alpha _i\lambda _i^kW_i$, and as 
$\lambda _i>0$ and $\lambda _i<\lambda _{i+1}$, the set of roots 
of $\Phi ^k(U)$ (resp. of $\Phi ^k(V)$) 
tends for $k\rightarrow \infty$ to the one of $W_j$ (resp. the one of 
$W_{j+1}$).~~~~~$\Box$\\ 

$4^0$. The numbers of positive roots of $\Phi ^k(U)$ and $\Phi ^k(V)$ are 
respectively $j$ and $j+1$. These roots interlace. This follows from 
Corollary~\ref{corcor} and from Theorem~\ref{preplit}.
The above proposition implies that  

\begin{equation}\label{star}
y_1\leq x_1\leq y_2\leq x_2\leq \cdots \leq y_j\leq x_j\leq y_{j+1}
\end{equation} 
The roots of $Q_j$ and the ones of $Q_{j+1}$ being distinct, there 
can be no two consecutive 
equalities in this string of inequalities. 

$5^0$. Suppose that there is an equality of the form $x_i=y_i$ 
or $x_i=y_{i+1}$. Denote by $x_0$ a root common for $Q_j$ and $Q_{j+1}$. 
The roots of $Q_j$ and $Q_{j+1}$ being simple one can find 
$(\theta _0,\mu _0)\in ({\bf R}^2\backslash \{ (0,0)\} )$ 
such that the polynomial 
$L:=\theta _0W_j+\mu _0 W_{j+1}$ has a multiple root at $x_0$. 

Consider the polynomial 
$M:=(\theta _0/\lambda _j)W_j+(\mu _0/\lambda _{j+1})W_{j+1}$. 
One has $L=\Phi (M)$. The polynomial $M$ can be considered as a limit 
of a sequence of polynomials of the form 
$M_s:=(\theta _0/\lambda _j)U_s+(\mu _0/\lambda _{j+1})W_{j+1}$, 
$s=1,2,\ldots$, 
where $U_s=x(x+1)^{n-2-j}T_{j,s}$ the roots of $T_{j,s}$ being positive, 
distinct and interlacing with the ones of $Q_{j+1}$. For each $s$ fixed 
the polynomial $M_s$ has at least $j$ positive distinct roots (and by 
Proposition~\ref{posroots} this is also the case of $L_s:=\Phi (M_s)$), 
hence $M$ and $L$ have each at least $j$ positive roots 
counted with multiplicity.

In the right-hand side of equality (\ref{factoriz}) with 
$P=M$ there are $n-j-3$ composition factors with $a_i=1$ which can be skipped. 
From the remaining ones there are $j$ or $j+1$ with $a_i<0$ and one or two 
with $a_i=0$. If there are (at least) 
two composition factors with the same $a_i<0$ 
(which is the case when $L=\Phi (M)$ has a multiple positive root), 
then their composition is 
a polynomial with all coefficients nonnegative, i.e. 
without change of the sign in the sequence 
of the coefficients. Hence the 
remaining composition factors with $a_i<0$ are not more than $j-1$ and their 
composition  
can be a polynomial with not more than $j-1$ changes of the sign 
in the sequence of the coefficients (see Proposition~\ref{compfactors}). 
The composition factor(s) $K_0$ add(s) 
no changes of sign in this sequence.

But then according to the Descartes rule the polynomial $M$ must have not more 
than $j-1$ positive roots counted with multiplicity 
which is a contradiction.~~~~~$\Box$

\section{Proof of Theorem~\protect\ref{tm2}\protect\label{prtm2}}

$1^0$. Recall that $W_{\nu}=x(x+1)^{n-2-\nu }Q_{\nu}$. 
Denote by $B$ a polynomial of the form 
$x(x+1)^{n-2-j}C$ where the polynomial $C$ is monic, degree $j$, 
self-reciprocal and has a root in each of the intervals 
$(z_i,z_{i+1})$ and $(z_{j+2-i},z_{j+3-i})$ for $i=1,\ldots ,[j/2]$. (If 
$j$ is odd, then $C(1)=0$.) 
In this sense we say that the roots of the polynomials $C$ and $Q_{j+2}$ 
satisfy the second interlacing property.

\begin{rem}\label{rrrr}
{\rm The polynomial $\Phi (B)$ is divisible by $x(x+1)^{n-2-j}$. This 
follows from part 1) of Proposition~\ref{propmult} 
and from Remark~\ref{remmult}. In the same 
way for any $k\in {\bf N}$ one can conclude that 
the polynomial $\Phi ^k(B)$ is divisible 
by $x(x+1)^{n-2-j}$.}
\end{rem}

\begin{lm}\label{generic}
One can choose the polynomial $C$ such that for any $k\in {\bf N}\cup \{ 0\}$ 
the polynomials $Q_{j+2}$ and $\Phi ^k(B)/(x(x+1)^{n-2-j})$ have no root 
different from $1$ in common.
\end{lm}

{\em Proof:}\\ 

Suppose that $j$ is even. Then $1$ is not a root of $C$ and not a root of 
$Q_{j+2}$. Fix $k\in {\bf N}$. The condition the polynomials 
$Q_{j+2}$ and $\Phi ^k(B)/(x(x+1)^{n-2-j})$ 
to have a root in common reads 
$\Delta _k:=$Res$(Q_{j+2},\Phi ^k(B)/(x(x+1)^{n-2-j}))=0$. 
This equality defines a proper algebraic subvariety in the space 
of the half of the coefficients of the polynomial $C$; half of them 
because $C$ is self-reciprocal. 

It is clear that $\Delta _k$ is a not 
identically constant polynomial in the half of the coefficients of $C$. 
Indeed, 

1) for $k=0$ this results from the fact that  
if one fixes all roots of $C$ except $x_0$ and $1/x_0$ the remaining ones 
being nonroots of $Q_{j+2}$, then $\Delta _0$ will equal $0$ precisely when 
$x_0$ and $1/x_0$ are roots of $Q_{j+2}$; 

2) for $k>0$ this follows from 1) and from  
the mapping $\Phi ^k$ being nondegenerate and preserving self-reciprocity, 
see \cite{Ko2}.


It is reasonable to consider the space $\tilde{R}$ 
of the roots of $C$ which belong to $(0,1)$ and not the space 
of all its roots. Indeed, the roots of $C$ belonging to $(1,\infty )$ 
equal $1/x_*$ where $x_*$ is a root in $(0,1)$; 
when deg$C$ is odd, then $C$ has a simple root at $(-1)$. 

The roots of $C$ are real and distinct. 
Therefore the following statement is true 
(see its proof at the end of the proof of the lemma): 

{\em Locally the mapping 
``~$U~:~$ roots of $C$ belonging to $(0,1)$ $\mapsto$ 
first half of its coefficients~'' is a diffeomorphism.}

Thus the condition 
$\Delta _k=0$ defines a proper algebraic subvariety in the space $\tilde{R}$.  
Its complement $Z_k$  
is a Zariski open dense 
subset of $\tilde{R}$. The intersection $\cap _{k=0}^{\infty}Z_k$ is nonempty. 
To choose $C$ as claimed by the lemma is the same as 
to choose the roots of $C$ which 
belong to $(0,1)$ from the set $\cap _{k=0}^{\infty}Z_k$.

If $j$ is odd, then $1$ is a simple root of $C$ and of $Q_{j+2}$ and 
one can apply the same reasoning to the polynomials 
$C/(x-1)$ and $Q_{j+2}/(x-1)$.

Prove the above statement. 
For each two roots of $C$ of the form $x_0$, $1/x_0$ 
consider the product $(x-x_0)(x-1/x_0)=x^2-sx+1$. Write down $C$ as a 
product of such degree $2$ polynomials. The respective quantities $s$ are all 
real and distinct. Therefore the mapping $Y$ which sends their tuple into the 
tuple of the values of their elementary symmetric polynomials is a 
diffeomorphism. The mapping $T$ sending these symmetric polynomials into the 
tuple of the first half of the coefficients of $C$ is affine 
nondegenerate triangular (to be checked directly). The mapping $V$ sending the 
roots of $C$ belonging to $(0,1)$ into the quantities $s$ 
is a diffeomorphism (easy to check). The mapping $U$ equals $T\circ Y\circ V$. 
Hence it is a local diffeomorphism.~~~~~$\Box$\\ 

$2^0$. From now till the end of the proof of Theorem~\ref{prtm2} 
we assume that the roots of $C$ satisfy the conclusion of Lemma~\ref{generic}. 

\begin{prop}\label{prop2int}
The polynomial $\Phi (B)$ is self-reciprocal. The roots of 
$\Phi (B)/(x(x+1)^{n-2-j})$ and 
$Q_{j+2}$ satisfy the second interlacing property.
\end{prop}

Before proving the proposition deduce Theorem~\ref{tm2} from it. For every 
$k\in {\bf N}$ the polynomial $\Phi ^k(B)$ is divisible by $x(x+1)^{n-2-j}$, 
see Remark~\ref{rrrr}. 
Applying $k$ times Proposition~\ref{prop2int}  
one shows that for any $k\in {\bf N}$ the polynomial 
$\Phi ^k(B)/(x(x+1)^{n-2-j})$ is self-reciprocal and the roots of 
$\Phi ^k(B)/(x(x+1)^{n-2-j})$ and $Q_{j+2}$ satisfy 
the second interlacing property. 

$3^0$. Set $Q_0:=1$. The polynomial $B$ can be presented 
in a unique way in the form 
$\sum _{\nu =0}^j\gamma _{\nu}W_{\nu}$, $\gamma _{\nu}\in {\bf R}$. Indeed, 
a priori it can be presented in the form 
$a(x+1)^{n-1}+\sum _{\nu =0}^{n-1}\gamma _{\nu}W_{\nu}$, see the 
eigenvectors of the mapping $\Phi$ before Remark~\ref{remlin}. 
It is divisible by $x$ and by $(x+1)^{n-2-j}$. Hence 
$a=\gamma _{j+1}=\gamma _{j+2}=\cdots =\gamma _{n-1}=0$. Moreover, 
$\gamma _j\neq 0$, otherwise $B$ must be divisible by $(x+1)^{n-1-j}$ which is 
false. 

One has $\Phi ^k(B)=\sum _{\nu =0}^j(\lambda _{\nu})^k\gamma _{\nu}W_{\nu}=
(\lambda _j)^k\sum _{\nu =0}^j(\lambda _{\nu}/\lambda _j)^k
\gamma _{\nu}W_{\nu}$. As $1\leq \lambda _{\nu}<\lambda _j$, for large values 
of $k$ the positive roots of the right-hand side are close to the ones 
of the polynomial 
$(\lambda _j)^k\gamma _jQ_j$. Passing to the limit when $k\rightarrow \infty$ 
one sees that the polynomial $Q_j$ has a root in each of the intervals 
$[z_i,z_{i+1}]$ and $[z_{j+2-i},z_{j+3-i}]$ for $i=1,\ldots ,[j/2]$.  

It is 
clear that for $j$ 
odd one has $Q_j(1)=0$. Therefore the second interlacing property will be 
proved if one manages to exclude the possibility $Q_j$ and $Q_{j+2}$ to have 
a common positive root different from $1$. We do this by analogy with 
$4^0$ and $5^0$ of the proof of Theorem~\ref{tm1}.

$4^0$. Suppose that such a root $x_0$ exists. 
Then $1/x_0$ is also such a root. 

Consider the polynomials $W_j$ and $W_{j+2}$ as degree $n$, not $n-1$ ones. 
They are self-reciprocal. One has simultaneously 
$x^{n}W_j(1/x)=W_j(x)$, $x^{n}W_{j+2}(1/x)=W_{j+2}(x)$ or 
$x^{n}W_j(1/x)=-W_j(x)$, $x^{n}W_{j+2}(1/x)=-W_{j+2}(x)$. Indeed, the 
parities of the multiplicities of their roots 
$1$ and $(-1)$ are the same. 
Therefore any linear combination $\theta W_j+\mu W_{j+2}$ is a self-reciprocal 
polynomial. 

The root $x_0$ (considered as a root of $W_j$ or $W_{j+2}$) is simple. 
Therefore there exists a couple 
$(\theta _0,\mu _0)\in ({\bf R}^2\backslash \{ (0,0)\} )$ 
such that the polynomial 
$W:=\theta _0W_j+\mu _0W_{j+2}$ has a multiple root at $x_0$. 
By self-reciprocity it 
has a multiple root at $1/x_0$ as well.   

$5^0$. The polynomial 
$H:=(\theta _0/\lambda _j)W_j+(\mu _0/\lambda _{j+1})W_{j+2}$ 
can be presented as a limit of 
a sequence of polynomials of the form 
$H_s:=(\theta _0/\lambda _j)G_s+(\mu _0/\lambda _{j+1})W_{j+2}$, 
$s=1,2,\ldots$, 
where $G_s=x(x+1)^{n-2-j}T_{j,s}$, deg$T_{j,s}=j$, $T_{j,s}$ 
is self-reciprocal, its roots are positive, 
distinct and interlacing with the ones of $Q_{j+2}$ in the sense of the second 
interlacing property. For each $s$ fixed 
the polynomial $H_s$ has at least $j$ positive distinct roots (and by 
Proposition~\ref{posroots} this is true also for the polynomial 
$W_s:=\Phi (H_s)$), hence $H$ and $W=\Phi (H)$ have each 
at least $j$ positive roots counted with multiplicity.

In the right-hand side of equality (\ref{factoriz}) with 
$P=H$ there are (except the $n-j-3$ composition factors $K_1$), $j$ or 
$j+1$ composition factors with $a_i<0$. There are among them two couples 
of equal factors. The composition of such a couple is a polynomial with all 
coefficients positive. Hence, it adds no change of the sign in the sequence of 
coefficients of the polynomial $W$. The remaining $\leq j-3$ composition 
factors with $a_i<0$ bring not more than $j-3$ such changes and the factor(s) 
$K_0$ bring(s) no change at all. By the Descartes rule the polynomial $H$ 
must have not more 
than $j-3$ positive roots counted with multiplicity 
which is a contradiction.~~~~~$\Box$\\

{\em Proof of Proposition~\ref{prop2int}:}\\ 

$1^0$. The polynomial $\Phi (B)$ is self-reciprocal because 
the matrix of the mapping $\Phi$ in the standard monomial basis is 
centre-symmetric, see Proposition~10 in \cite{Ko2}. 

$2^0$. Consider any linear combination 
$\Xi :=\theta (x+1)^2C+\mu Q_{j+2}$, 
$(\theta ,\mu )\in ({\bf R}^2\backslash \{ (0,0)\} )$. This degree $j+2$ 
polynomial has $j+2$ (resp. $j$) positive distinct roots for $\theta =0$ 
(resp. for $\mu =0$). For $\theta \neq 0\neq \mu$ it changes sign at any two 
consecutive positive roots of $Q_{j+2}$ 
which are both smaller or both greater than $1$. 
If $j$ is odd, then $\Xi (1)=0$. 
Hence for all $(\theta ,\mu )$, the polynomial $\Xi$ has at least $j$ 
real positive distinct roots (and at most one complex conjugate couple). 

$3^0$. Hence the polynomial 

\[ \tilde{B}:=\Phi (x(x+1)^{n-j-4}\Xi )=
\Phi (\theta B+\mu W_{j+2})=\theta \Phi (B)+\mu \lambda _{j+4}W_{j+2}\]  
has a root at $0$ (part 1) of Proposition~\ref{propmult}), 
at least $j$ distinct 
positive roots (Proposition~\ref{posroots}), 
a root at $(-1)$ of multiplicity at least $n-j-4$ (Remark~\ref{remmult}). The 
remaining two roots are either complex or real; in the latter case one or 
both of them can coincide with some of the other roots.  
In particular, they can equal $(-1)$.

$4^0$. Consider the rational function $g:=\Phi (B)/W_{j+2}$. 
Lemma~\ref{generic} allows 
to suppose that 
$\Phi (B)$ and $W_{j+2}$ have no positive root different from $1$ in common. 
Then $g(x)$ has no 
critical point for $x\in (0,1)\cup (1,\infty )$. Indeed, self-reciprocity 
of $\Phi (B)$ and $W_{j+2}$ would imply that there 
are in fact two different such critical points, at $x_0$ and $1/x_0$, 
with equal critical 
values. If they are the only ones at their level set and if they are 
Morse ones, then  
for two different values of $(\theta ,\mu )$ the number of distinct 
positive roots of the polynomial $\tilde{B}$ 
changes by $4$ which by $2^0$ is impossible because 
the real positive roots of $\tilde{B}$ are not more than $j+2$ 
and not less than $j$. 

If the critical points are either more than two on their level set 
or they are not Morse ones, 
then for some $(\theta ,\mu )$ the number of distinct positive roots 
of the polynomial $\tilde{B}$ 
is smaller than $j$ -- a contradiction with $2^0$ again.

$5^0$. Therefore $g$ must be monotonous between each two consecutive 
positive roots of $Q_{j+2}$ which are smaller than $1$, and between 
each two consecutive positive roots of $Q_{j+2}$ which are greater than $1$. 
The function $g$ has a single critical point, at $1$, which is a Morse one. 
This implies that the roots of $\Phi (B)/(x(x+1)^{n-2-j})$ 
(which are the level set 
$\{ g=0\}$ without the double root at $(-1)$) and the ones of $Q_{j+2}$ 
satisfy the second interlacing property.~~~~~$\Box$


\section{On Gegenbauer and Narayana polynomials\protect\label{Geg}}

In this paper we define {\em the Gegenbauer polynomial} $G_n$ of degree $n$ 
($n\geq 3$) 
as a polynomial in one variable which is divisible by its 
second derivative and whose first three coefficients equal $1$, $0$, $-1$. 
It is easy to show that these conditions define a unique polynomial which is 
strictly hyperbolic and which is even (resp. odd) 
when $n$ is even (resp. odd). Hence its nonconstant 
derivatives which are polynomials of even (resp. odd) degree are even 
(resp. odd) polynomials. The general definition of Gegenbauer polynomials 
$C_n^{(\lambda )}$ is different from the above one and depends on a  
parameter $\lambda$. Our definition corresponds to the particular case 
$\lambda =-1/2$. (For more details about Gegenbauer polynomials 
see Chapter 22 in \cite{AbSt}.) 
Up to a rescaling and a nonzero constant factor the 
polynomial $G_n'$ is the Legendre polynomial of degree $n-1$.

Denote by $\zeta _1<\cdots <\zeta _l<0$ 
the negative roots of $G_n^{(k)}$ where $l=[(n-k)/2]$. Denote by 
$\mu _1<\cdots <\mu _{l-1}<0$ the negative roots of $G_n^{(k+2)}$. 

\begin{tm}\label{Gegtm}
One has $\mu _i\in (\zeta _i,\zeta _{i+1})$, $i=1,\ldots ,l-1$.
\end{tm}

{\em Proof:}\\ 

The polynomial $G_n$ satisfies the condition $G_n=(x^2-a^2)G_n''$ where 
$\pm a$ are its roots with greatest absolute values. Differentiate this 
equality $k$ times using the Leibniz rule:

\[ G_n^{(k)}=(x^2-a^2)G_n^{(k+2)}+2kxG_n^{(k+1)}+k(k-1)G_n^{(k)}~.\] 
Hence 
$(\zeta _i^2-a^2)G_n^{(k+2)}(\zeta _i)+2k\zeta _iG_n^{(k+1)}(\zeta _i)=0$. 
Recall that $\zeta _i^2-a^2<0$, $\zeta _i<0$. This means that the signs 
of $G_n^{(k+2)}(\zeta _i)$ and $G_n^{(k+1)}(\zeta _i)$ are opposite from 
where the theorem follows easily.~~~~~$\Box$\\  

For Narayana polynomials the following recurrence relation holds 
(see \cite{Su}):

\begin{equation}\label{Narayanarecur}
(n+1)N_n(x)=(2n-1)(1+x)N_{n-1}(x)-(n-2)(x-1)^2N_{n-2}(x)
\end{equation}




\begin{rem}\label{Nrem}
{\rm It is shown in paper \cite{KoMaSh} (see part 1) of Corollary~7 there) 
that the roots of the Narayana polynomials $N_n$ and $N_{n-1}$ interlace, i.e. 
satisfy the first interlacing property.}
\end{rem} 

\begin{tm}
The zeros of the Narayana polynomials $N_n$ and $N_{n-2}$ satisfy the second 
interlacing property.
\end{tm}

Indeed, it suffices to show this for the roots greater than $1$, 
for the others it will follow from the self-reciprocity of Narayana 
polynomials. The claimed property results from equation 
(\ref{Narayanarecur}) -- if $\xi >1$ is the greatest 
root of $N_{n-2}$, then $N_{n-1}(\xi )<0$ (see Remark~\ref{Nrem}). By 
equality (\ref{Narayanarecur}) the signs of 
$N_n(\xi )$ and $N_{n-1}(\xi )$ are the same. As $N_{n-1}(\xi )<0$, 
one has $N_n(\xi )<0$, i.e. $\xi$ is greater than the 
greatest but one root of $N_n$. In the same way the second 
interlacing property is proved for the other roots 
of $N_{n-2}$ which are $>1$.~~~~~$\Box$

Author's address: Universit\'e de Nice, Laboratoire de Math\'ematiques, 
Parc Valrose, 06108 Nice Cedex 2, France

e-mail: kostov@math.unice.fr


\begin{thebibliography}{Dillo 83}
\bibitem{Ko3} {\sc Kostov V.P.} Teorema realizatsii v kontekste 
kompozitsii Shura-Sege 
(A realization theorem in the context of the Schur-Szeg\H{o} composition), 
Functional Analysis and its Applications (accepted).
\end{thebibliography}
\end{document}